\documentclass[titlepage,11pt]{article}
\usepackage{latexsym}
\usepackage{amsfonts}
\usepackage{amsmath}
\newcommand\blackslug{\hbox{\hskip 1pt \vrule width 4pt height 8pt depth 1.5pt
        \hskip 1pt}}
\newcommand\bbox{\hfill \quad \blackslug \medbreak}

\def\l{,\ldots,}
\def\Reals{{\bf R_+}}
\def\mymod{\text{ mod }}

\title{Cycles in dense digraphs}
\author{Maria Chudnovsky\thanks{This research was conducted while the author served
as a Clay Mathematics Institute Research Fellow.}\\ 
Paul Seymour\thanks{Supported by ONR grant N00014-04-1-0062, and NSF grant
DMS03-54465.}\\
Blair Sullivan\thanks{This research was performed under an appointment to the Department of
Homeland Security (DHS) Scholarship and Fellowship Program.}\\
Princeton University, Princeton NJ 08544}

\date{July 24, 2006; revised \today}

\newtheorem{thm}{}[section]

\newcommand{\Proof}{\noindent{\bf Proof.}\ \ }

\begin{document}
\maketitle
\begin{abstract}
Let $G$ be a digraph (without parallel edges) such that every directed cycle has length at
least four; let $\beta(G)$ denote the size of the smallest subset $X\subseteq E(G)$
such that $G\setminus X$ has no directed cycles, and let $\gamma(G)$ be the number of unordered pairs $\{u,v\}$ 
of vertices such that $u,v$ are nonadjacent in $G$. It is easy to see that if $\gamma(G) = 0$
then $\beta(G) = 0$; what can we say about $\beta(G)$ if $\gamma(G)$ is bounded? 

We prove that
in general $\beta(G)\le \gamma(G)$. We conjecture that in fact $\beta(G)\le \frac12 \gamma(G)$ 
(this would be best possible if true),
and prove this conjecture in two special cases: 
\begin{itemize}
\item when $V(G)$ is the union of two cliques,
\item when the vertices of $G$ can be arranged in a circle
such that if distinct $u,v,w$ are in clockwise order and $uw$ is a (directed) edge, then so are both
$uv,vw$.
\end{itemize}
\end{abstract}

\section{Introduction}
We begin with some terminology. 
All digraphs in this paper are finite and have no parallel edges; and for a digraph $G$, $V(G)$ and $E(G)$
denote its vertex- and edge-sets. The members of $E(G)$ are ordered pairs of vertices, and we abbreviate
$(u,v)$ by $uv$. For integer $k\ge 0$, let us say a digraph $G$ is {\em $k$-free} if there is no directed cycle of $G$
with length at most $k$. 
A digraph is {\em acyclic} if it has no directed cycle.

We are concerned here with $3$-free digraphs.
It is easy to see that every $3$-free tournament is acyclic, and one might hope that 
every $3$-free digraph that is ``almost'' a tournament is ``almost'' acyclic. That is the topic of this paper.

More exactly, for a digraph $G$, let $\gamma(G)$ be the number of unordered pairs $\{u,v\}$ of distinct vertices $u,v$ that are
nonadjacent in $G$ (that is, both $uv,vu\notin E(G)$). Thus, every $2$-free digraph $G$
can be obtained from a tournament by deleting $\gamma(G)$ edges. 
Let $\beta(G)$ denote the minimum cardinality of a set $X\subseteq E(G)$ such that $G\setminus X$ is acyclic.
We already observed that every $3$-free digraph with $\gamma(G) = 0$ satisfies $\beta(G) = 0$, and our first result
is an extension of this.

\begin{thm}\label{betathm}
If $G$ is a $3$-free digraph then $\beta(G)\le \gamma(G)$.
\end{thm}
\Proof
We proceed by induction on $|V(G)|$, and we may assume that $V(G)\ne \emptyset$. 
Let us say a {\em $2$-path} is a triple $(x,y,z)$ such that $x,y,z\in V(G)$
are distinct, and $xy,yz\in E(G)$, and $x,z$ are nonadjacent. For each vertex $v$, let $f(v)$ denote the number of   
$2$-paths $(x,y,z)$ with $x = v$, and let $g(v)$ be the number of $2$-paths $(x,y,z)$ with $y = v$. Since $V(G)\ne \emptyset$
and $\sum_{v\in V(G)}f(v) = \sum_{v\in V(G)}g(v)$, there exists $v\in V(G)$ such that $f(v)\le g(v)$. 
Choose some such vertex $v$, and let $A,B,C$
be the set of all vertices $u\ne v$ such that $vu\in E(G)$, $uv\in E(G)$, and $uv,vu\notin E(G)$ respectively.
Thus the four sets $A,B,C,\{v\}$ are pairwise disjoint and have union $V(G)$.
Let $G_1,G_2$ be the subdigraphs of $G$ induced on $A$ and on $B\cup C$ respectively.
Since $g(v)$ is the number of pairs $(a,b)$ with $a\in A$ and
$b\in B$ such that $a,b$ are nonadjacent, it follows that $\gamma(G)\ge \gamma(G_1)+\gamma(G_2) + g(v)$.
From the inductive hypothesis, 
$\beta(G_1)\le \gamma(G_1)$ and $\beta(G_2)\le \gamma(G_2)$; for $i = 1,2$, choose $X_i\subseteq E(G_i)$
with $|X_i|\le \beta(G_i)$ such that $G_i\setminus X_i$ is acyclic. 
Let $X_3$ be the set of all edges $ac\in E(G)$ with $a\in A$ and $c\in C$; thus $|X_3| = f(v)$.
Since there is no edge $xy\in E(G)$ with $x\in A$ and $y\in B$ (because $G$ is $3$-free), it follows that
every edge $xy$ with $x\in A$ and $y\in \{v\}\cup B\cup C$ belongs to $X_3$, and so $G\setminus X$ is acyclic, where
$X = X_1\cup X_2\cup X_3$. Hence
$$\beta(G)\le |X| = |X_1|+|X_2|+|X_3| = \beta(G_1)+\beta(G_2)+f(v)\le \gamma(G_1)+\gamma(G_2)+g(v)\le \gamma(G).$$
This proves \ref{betathm}.~\bbox

Unfortunately, \ref{betathm} does not seem to be sharp, and we believe that the following holds.

\begin{thm}\label{mainconj}
{\bf Conjecture.} If $G$ is a $3$-free digraph then $\beta(G)\le \frac12 \gamma(G)$.
\end{thm}

If true, this is best possible for infinitely many values of $\gamma(G)$. For instance, let $G$
be the digraph with vertex set $\{v_1\l v_{4n}\}$, and with edge set as follows (reading subscripts modulo $4n$):
\begin{itemize}
\item $v_iv_j\in E(G)$ for all $i,j,k$ with $1\le k\le 4$ and $(k-1)n<i<j\le kn$
\item $v_iv_j\in E(G)$ for all $i,j,k$ with $1\le k\le 4$ and $(k-1)n<i\le kn<j\le (k+1)n$.
\end{itemize}
It is easy to see that this digraph $G$ is $3$-free, and satisfies $\beta(G) = n^2$ (certainly $\beta(G)\ge n^2$
since $G$ has $n^2$ directed cycles that are pairwise edge-disjoint), and $\gamma(G) = 2n^2$. 

The reason for our interest in \ref{mainconj} was originally its application to the Caccetta-H\"{a}ggkvist
conjecture~\cite{CH}.
A special case of that conjecture asserts the following:
\begin{thm}\label{CHconj}
{\bf Conjecture.} If $G$ is a $3$-free digraph with $n$ vertices, then some vertex
has outdegree less than $n/3$.
\end{thm}
This is a challenging open question and has received a great deal of attention.
Any counterexample to \ref{CHconj} satisfies $\gamma(G)\le \frac12 |E(G)|$, so our conjecture
\ref{mainconj} would tell us that $\beta(G)\le \frac14 |E(G)|$, and this would perhaps be useful information
towards solving \ref{CHconj}.

We have not been able to prove \ref{mainconj} in general, and in this paper we prove two partial results,
that \ref{mainconj} holds for every $3$-free digraph $G$ such that either
\begin{itemize}
\item $V(G)$ is the union of two cliques, or
\item the vertices of $G$ can be arranged in a circle
such that if distinct $u,v,w$ are in clockwise order and $uw\in E(G)$, then 
$uv,vw\in E(G)$.
\end{itemize}

The first result is proved in \ref{twocliques}, and the second in \ref{circintthm}.
Incidentally, Kostochka and Stiebitz~\cite{KS} proved 
that in any minimal counterexample to \ref{mainconj}, every vertex is nonadjacent to at least three other vertices,
and the conjecture is true for all digraphs with at most $8$ vertices.

\section{A distant relative of the four functions theorem}

In this section we prove a result that we apply in the next section.
We begin with an elementary lemma. ($\Reals$ denotes the set of nonnegative real numbers.)
\begin{thm}\label{quad}
If $a_1,a_2,c_1,c_2,d_1,d_2\in \Reals$ and $a_k^2\le c_kd_k$ for $k = 1,2$, then $(a_1+a_2)^2\le (c_1+d_1)(c_2+d_2)$.
\end{thm}
\Proof
If say $c_1 = 0$, then since $a_1^2\le c_1d_1$, it follows that $a_1 = 0$, and so
$$(a_1+a_2)^2 = a_2^2\le c_2d_2\le (c_1+c_2)(d_1+d_2)$$
as required. We may therefore assume that $c_1,c_2$ are nonzero. Now
\begin{eqnarray*}
	(c_1+c_2)(d_1+d_2)
 		&=& c_1d_1+c_1d_2+c_2d_1+c_2d_2\\
                &\ge& a_1^2 + c_1(a_2^2/c_2) + c_2(a_1^2/c_1) + a_2^2\\
                &=& (a_1+a_2)^2 + c_1c_2(a_2/c_2 - a_1/c_1)^2 \\
                &\ge& (a_1+a_2)^2.
\end{eqnarray*}
This proves \ref{quad}.~\bbox

Before the main result of this section we must set up some notation.
Let $m,n\ge 1$ be integers, and let $V$ denote the set of all pairs $(i,j)$ of integers with 
$1\le i\le m$ and $1\le j\le n$. If $f:V\rightarrow\Reals$,
and $X\subseteq V$, we define
$f(X)$ to mean $\sum_{x\in X}f(x)$. For $(i,j),(i',j')\in V$, we say that $(i',j')$ {\em dominates} $(i,j)$
if $i<i'$ and $j<j'$. 
Let $a,b:V\rightarrow \Reals$ be functions. We say that $b$ {\em dominates} $a$ if 
\begin{itemize}
\item $a(V) = b(V)$
\item for all $X,Y\subseteq V$, if $a(X) + b(Y)> a(V)$ then there exist $x\in X$ and $y\in Y$ such that
$y$ dominates $x$.
\end{itemize}

The main result of this section is the following. 
(It is reminiscent of the ``four functions'' theorem of Ahlswede and Daykin~\cite{AD},
but we were not able to derive it from that theorem.)

\begin{thm}\label{ourfunctions} Let $m,n\ge 1$ be integers, let $V$ be as above, and let $a,b,c,d$ be functions from 
$V$ to $\Reals$,
satisfying the following:
\begin{enumerate}
\item $a(i,j)b(i',j')\le c(i',j)d(i,j')$ for $1\le i<i'\le m$ and $1\le j<j'\le n$, and
\item $b$ dominates $a$.
\end{enumerate}
Then $a(V)b(V)\le c(V)d(V)$.
\end{thm}
\Proof
We proceed by induction on $m+n$. Let $Q$ be the set of all quadruples $(a,b,c,d)$ 
of functions from $V$ to $\Reals$ that satisfy conditions 1 and 2 above.
We say that $(a,b,c,d)\in Q$ is {\em good} if 
$$a(V)b(V)\le c(V)d(V).$$
Thus, we need to show that every member of $Q$ is good.
Certainly if $m = 1$ or $n = 1$ then condition 2 
implies that $a(V) = b(V) = 0$,
and therefore $(a,b,c,d)$ is good; so we may assume that $m,n\ge 2$.
\\
\\
(1) {\em If $(a,b,c,d)\in Q$ then $b(i,1) = 0$ for $1\le i\le m$, and $a(m,j) = 0$ for $1\le j\le n$.}
\\
\\
For let $X = V$, and let $Y$ be the set of all pairs $(i,1)$ with $1\le i\le m$. There do not
exist $x\in X$ and $y\in Y$ such that $y$ dominates $x$, and since $b$ dominates $a$ it follows that
$a(X) + b(Y)\le a(V)$. Since $a(X) = a(V)$ we deduce that $b(Y) = 0$. This proves the first statement, and the
second follows similarly. This proves (1).
\\
\\
(2) {\em If $(a,b,c,d)\in Q$ and $a(i,1) = 0$ for all $i\in \{1\l m\}$ then $(a,b,c,d)$ is good.}
\\
\\
This follows from (1) and the inductive hypothesis applied to the restriction of $a,b,c,d$ to the set of all $(i,j)\in V$
with $j>1$ (relabeling appropriately).
\bigskip

For  $(a,b,c,d)\in Q$, let us define its {\em margin} to be the number of pairs $(i,j)$ such that 
either $j= 1$ and $a(i,j)>0$, or
$i = m$ and $b(i,j)>0$. For fixed $m,n$ we proceed by induction on the margin. Thus, we assume that
$t\ge 0$ is an integer, and every $(a,b,c,d)\in Q$ with margin smaller than $t$ is good.
We must show that every $(a,b,c,d)\in Q$ with margin $t$ is good.
\\
\\
(3) {\em Let $(a,b,c,d)\in Q$ with margin $t$, and suppose that there exist $X,Y\subseteq V$ such that
\begin{itemize}
\item $a(X) + b(Y) = a(V)$
\item there do not exist $x\in X$ and $y\in Y$ such that $y$ dominates $x$
\item there exists $i\in \{1\l m\}$ such that $(i,1)\notin X$ and $a(i,1)>0$, and there exists
$j\in \{1\l n\}$ such that $(m,j)\notin Y$ and $b(m,j)>0$.
\end{itemize}
Then $(a,b,c,d)$ is good.}
\\
\\
Let $A_1 = X$ and $A_2 = V\setminus X$. Let $B_1 = V\setminus Y$ and $B_2 = Y$.
For $k = 1,2$, let $C_k$ be the set of all pairs $(i',j)\in V$ such that there exist $i,j'$ with $i<i'$ and $j<j'$
and $(i,j)\in A_k$ and $(i',j')\in B_k$; and let $D_k$ be the set of all pairs $(i,j')$ 
such that there exist $i',j$ with $i<i'$ and $j<j'$
and $(i,j)\in A_k$ and $(i',j')\in B_k$. We observe first that $C_1\cap C_2 = \emptyset$; for suppose that
$(i',j)\in C_1\cap C_2$. Since $(i',j)\in C_1$, there exists $i<i'$ such that $(i,j)\in X$; and since
$(i',j)\in C_2$, there exists $j'>j$ such that $(i',j')\in Y$. But then $(i',j')\in Y$ dominates $(i,j)\in X$, contradicting
the second hypothesis about $X,Y$. This proves that $C_1\cap C_2 = \emptyset$, and similarly $D_1\cap D_2 = \emptyset$.
For $k = 1,2$, and $x\in V$, define $a_k(x) = a(x)$ if $x\in A_k$, and $a_k(x) = 0$ otherwise.
Define $b_k(x),c_k(x),d_k(x)$ similarly. 
Since $a_1(V)+a_2(V), b_1(V)+b_2(V)$ and $a_1(V)+b_2(V)$ all equal $a(V)$, it follows that $a_1(V) = b_1(V)$ and
$a_2(V) = b_2(V)$.
We claim that $(a_k,b_k,c_k,d_k)\in Q$ for $k = 1,2$.
To see this, let $i<i'$ and $j<j'$; we must show first that
$a_k(i,j)b_k(i',j')\le c_k(i',j)d_k(i,j')$.
Hence we may assume that $a_k(i,j)$ and $b_k(i',j')\ne 0$, and therefore $(i,j)\in A_k$ and $(i',j')\in B_k$.
From the definition of $C_k,D_k$ it follows that $(i',j)\in C_k$ and $(i,j')\in D_k$. Hence
$a_k(i,j) = a(i,j)$, and $b_k(i',j') = b(i',j')$, and $c_k(i',j) = c(i',j)$, and $d_k(i,j') = d(i,j')$; and since
$a(i,j)b(i',j')\le c(i',j)d(i,j')$, this proves the claim. Second, we must show that $b_k$ dominates $a_k$.
We have already seen that $a_k(V) = b_k(V)$. Let $X',Y'\subseteq V$
with $a_k(X') + b_k(Y')> a_k(V)$; we must show that
there exist $x\in X'$ and $y\in Y'$ such that $y$ dominates $x$. From the symmetry we may assume that $k = 1$.
Now $a(X\cap X') = a_k(X')$, and $b(Y\cup Y') = b(Y) + b_k(Y')$, and so
$$a(X\cap X') + b(Y\cup Y') = a_k(X') +  b(Y) + b_k(Y')> a_k(V) + b(Y) = a(X) + b(Y) = a(V).$$
Since $b$ dominates $a$, there exist $x\in X\cap X'$ and $y\in Y\cup Y'$ such that $y$ dominates $x$. No vertex in
$Y$ dominates a vertex in $X$, from the choice of $X,Y$, and it follows that $y\in Y'$,
as required. This proves that $b_k$ dominates $a_k$, and consequently $(a_k,b_k,c_k,d_k)\in Q$, for $k = 1,2$.

We claim that for $k = 1,2$, the margin of $(a_k,b_k,c_k,d_k)$ is less than $t$. For from the third
hypothesis about $X,Y$,
there exists $i\in \{1\l m\}$ such that $a(i,1)>0$ and $(i,1)\notin X$ (and hence $a_1(i,1) = 0$); this
shows that the margin of $(a_1,b_1,c_1,d_1)$ is less than that of $(a,b,c,d)$, and so less than $t$. Also, there exists
$j\in \{1\l n\}$ such that $b(m,j)>0$ and $(m,j)\notin Y$; and so similarly the margin of
$(a_2,b_2,c_2,d_2)$ is less than $t$. Hence from the second inductive hypothesis, we deduce that
$a_k(V)b_k(V)\le c_k(V)d_k(V)$ for $k = 1,2$.
But $a_k(V) = b_k(V)$ for $k = 1,2$; thus $a_k(V)^2\le c_k(V)d_k(V)$ for $k = 1,2$.
Since $a_1(V)+a_2(V) = a(V) = b(V)$ and since $c(V) \ge c_1(V)+c_2(V)$ (because $C_1\cap C_2 = \emptyset$),
and similarly $d(V) \ge d_1(V)+d_2(V)$, it suffices to show that
$$(a_1(V)+a_2(V))^2\le (c_1(V)+c_2(V))(d_1(V)+d_2(V)),$$
and this follows from \ref{quad}.
This proves (3).
\\
\\
(4) {\em If $(a,b,c,d)\in Q$ with margin $t$, and there exists $j\ge 3$ such that $b(m,j)>0$, then 
$(a,b,c,d)$ is good.}
\\
\\
For let $\epsilon$ satisfy $0\le \epsilon\le 1$. For $1\le i\le m$, define 
\begin{eqnarray*}
a_1(i,1) &=& (1-\epsilon)a(i,1)\\
a_1(i,2) &=& \epsilon a(i,1) + a(i,2)\\
a_1(i,j) &=& a(i,j) \text{ for } 3\le j\le n\\
c_1(i,1) &=& (1-\epsilon)c(i,1)\\
c_1(i,2) &=& \epsilon c(i,1) + c(i,2)\\
c_1(i,j) &=& c(i,j) \text{ for } 3\le j\le n.
\end{eqnarray*}

Since $b$ dominates $a$, by compactness we may choose $\epsilon\le 1$ maximum such that $b$ dominates $a_1$.
We claim that $(a_1,b,c_1,d)\in Q$; for let $i<i'$ and $j<j'$. We must check that
$a_1(i,j)b(i',j')\le c_1(i',j)d(i,j')$.
If $j = 1$, then
$$a_1(i,j)b(i',j') = (1-\epsilon)a(i,1)b(i',j')$$
and
$$c_1(i',j)d(i,j') = (1-\epsilon)c(i,1)d(i,j'),$$
and since $a(i,1)b(i',j')\le c(i,1)d(i,j')$ it follows that $a_1(i,j)b(i',j')\le c_1(i',j)d(i,j')$
as required.
If $j = 2$, then 
$$a_1(i,j)b(i',j') = (\epsilon a(i,1) + a(i,2))b(i',j')$$
and 
$$c_1(i',j)d(i,j') = (\epsilon c(i,1) + c(i,2))d(i,j'),$$
and since $a(i,1)b(i',j')\le c(i,1)d(i',j')$ and $a(i,2)b(i',j')\le c(i,2)d(i',j')$, it follows
that $a_1(i,j)b(i',j')\le c_1(i',j)d(i,j')$ as required. 
Finally, if $j>2$ the claim is clear, since $a_1(i,j) = a(i,j)$ and $c_1(i',j) = c(i',j)$. 
This proves that $(a_1,b,c_1,d)\in Q$.

We claim that $(a_1,b,c_1,d)$ is good.
If $\epsilon = 1$, then $a_1(i,1) = 0$ for $1\le j\le m$, and therefore 
$(a_1,b,c_1,d)$ is good by (2). We may therefore assume that $\epsilon<1$. From the maximality of $\epsilon$,
there exist $X,Y\subseteq V$ such that 
\begin{itemize}
\item there does not exist $x\in X$ and $y\in Y$ such that $y$ dominates $x$
\item $a_1(X)+b(Y) = a_1(V)$
\item for some $i$ with $1\le i\le m$, $(i,1)\notin X$ and and $(i,2)\in X$ and $a(i,1)>0$.
\end{itemize}
(The third statement follows from the fact that increasing $\epsilon$ will cause $a_1(X)$ strictly to increase.)
Now we recall that there exists $j\ge 3$ such that $b(m,j)>0$. Since $(i,2)\in X$ is dominated by $(m,j)$ (for
$i<m$ by (1), since $a(i,2)>0$), it follows that $(m,j)\notin Y$. But then $(a_1,b,c_1,d)$ satisfies
the hypotheses of (3), and therefore $(a_1,b,c_1,d)$ is good. This proves the claim.

Since $a_1(V) = a(V) $ and $c_1(V) = c(V) $, we deduce that $(a,b,c,d)$ is good. This proves (4).
\bigskip

Now let $(a,b,c,d)\in Q$ with margin $t$; we shall prove that it is good. By (4) we may assume that
$b(m,j) = 0$ for $3\le j\le m$, and similarly that $a(i,1) = 0$ for $1\le i\le m-2$. Since $a(m,1) = b(m,1) = 0$
by (1), it follows that $a(i,1)\ne 0$ only if $i = m-1$, and $b(m,j)\ne 0$ only if $j = 2$. Let $X = \{(m-1,1)\}$
and let $Y$ be the set of all $(i,j)\in V$ with $i<m$; then there do not exist $x\in X$ and $y\in Y$
such that $y$ dominates $x$. Consequently $a(X) + b(Y)\le a(V)$. But $a(X) = a(m-1,1)$ and $b(Y)\ge a(V)-b(m,2)$, and so
$a(m-1,1)\le b(m,2)$. Similarly the reverse inequality holds, and so $a(m-1,1)= b(m,2)$. For $(i,j)\in V$, if
either $i = m$ or $j = 1$, define 
$$a_1(i,j) = b_1(i,j) = c_1(i,j) = d_1(i,j) = 0.$$
If $i<m$ and $j>1$ let $a_1(i,j) = a(i,j),b_1(i,j) = b(i,j)$, and $c_1(i,j) = c(i,j)$; and let $d_1(i,j) = d(i,j)$
except that $d_1(m-1,2) = 0$. We claim that $(a_1,b_1,c_1,d_1)\in Q$. For let $i<i'$ and $j<j'$. We must check that
$a_1(i,j)b_1(i',j')\le c_1(i',j)d_1(i,j')$.
If $i'<m$ and $j>1$ then $a_1(i,j) = a(i,j)$ and so on, and the claim is clear. If $i' = m$ or $j = 1$ 
then $a_1(i,j)b_1(i',j') = 0$ and again the claim is clear. Thus $a_1(i,j)b_1(i',j')\le c_1(i',j)d_1(i,j')$.
Next we must check that $b_1$ dominates $a_1$. Certainly
$$a_1(V) = a(V) - a(m-1,1) = b(V) - b(m,2) = b_1(V).$$
Let $X,Y\subseteq V$ such that $a_1(X) + b_1(Y)>a_1(V)$. We must show that there exist $x\in X$ and $y\in Y$
such that $y$ dominates $x$. We may therefore assume that $a_1(x)>0$ for all $x\in X$, and $b_1(y)>0$
for all $y\in Y$. In particular, since $a_1(m-1,1) = b_1(m,2) = 0$, it follows that $(m-1,1)\notin X$ and $(m,2)\notin Y$.
Let $X' = X\cup \{(m-1,1)\}$. Then $a(X') = a_1(X)+a(m-1,1)$, and so
$$a(X') + b(Y) = a_1(X)+a(m-1,1) + b(Y) > a_1(V) + a(m-1,1)= a(V).$$
Hence there exist $x\in X'$ and $y\in Y$ such that $y$ dominates $x$. If $x = (m-1,1)$, then $y = (m,j)$ for some $j>1$,
and therefore $b_1(y) = 0$, a contradiction, since $b_1(y)>0$ for all $y\in Y$. Thus $x\ne (m-1,1)$, and so $x\in X$,
as required. This proves that $b_1$ dominates $a_1$. 

By (2), $(a_1,b_1,c_1,d_1)$ is good, and so 
$a_1(V)b_1(V)\le c_1(V)d_1(V)$. Moreover,
$$a(m-1,1)b(m,2)\le c(m,1)d(m-1,2)$$
and $b(m,2) = a(m-1,1)$, and so $a(m-1,1)^2\le c(m,1)d(m-1,2)$.
Hence \ref{quad} implies that
$$(a_1(V) + a(m-1,1))^2\le (c_1(V) +c(m,1))(d_1(V)+d(m-1,2)).$$
But $a(V) = a_1(V) + a(m-1,1) = b(V)$, and $c(V) \ge c_1(V) +c(m,1)$, and $d(V)\ge d_1(V)+d(m-1,2)$;
and it follows that $(a,b,c,d)$ is good. This completes the inductive proof that every member of $Q$ is good,
and so proves \ref{ourfunctions}.~\bbox

\section{The two cliques result}

In this section we prove the following.
\begin{thm}\label{twocliques}
Let $G$ be a $3$-free digraph and let $M,N$ be a partition of $V(G)$ such that $M,N$ are both 
cliques of $G$. Then there is a set $X\subseteq E(G)$ such that every member of $X$ has one end in $M$ and one
end in $N$, and $|X|\le \frac12 \gamma(G)$, and $G\setminus X$ is acyclic. In particular, $\beta(G)\le \frac12 \gamma(G)$.
\end{thm}
\Proof The second assertion follows immediately from the first, so we just prove the first.
Since the restriction of $G$ to $M$ is a $3$-free tournament, we can number $M= \{u_1\l u_m\}$ such that
$u_iu_{i'}\in E(G)$ for $1\le i<i'\le m$. The same holds for $N$, but it is convenient to number its members in reverse order;
thus we assume that $N = \{v_1\l v_n\}$, where $v_{j'}v_{j}\in E(G)$ for $1\le j<j'\le n$. 
Let $V$ be the set of all
pairs $(i,j)$ with $1\le i\le m$ and $1\le j\le n$.
For $a = (i,j)\in V$ and $b = (i',j')\in V$, let us say that $(a,b)$ is a {\em cross} if 
$v_ju_i, u_{i'}v_{j'}\in E(G)$ and $1\le i<i'\le m$ and $1\le j<j'\le n$.
Let $A_0$ be the set of all edges of $G$ from $N$ to $M$, and $B_0$ the set of all edges from $M$ to $N$.
Let $k$ be the minimum cardinality of a subset $X\subseteq A_0\cup B_0$ such that $G\setminus X$ is acyclic.
(Such a number exists since $G\setminus (A_0\cup B_0)$ is acyclic.)
\\
\\
(1) {\em There are $k$ crosses $(a_1,b_1)\l (a_k,b_k)$ such that $a_1\l a_k,b_1\l b_k$ are all distinct.}
\\
\\
For suppose not. 
Let $H$ be the bipartite graph with vertex set $A_0\cup B_0$, in which $v_ju_i\in A_0$ and $u_{i'}v_{j'}\in B_0$ are adjacent
if $((i,j),(i',j'))$ is a cross. Then $H$ has no $k$-edge matching, and so by K\"onig's theorem, there 
exists $X\subseteq A_0\cup B_0$
with $|X|<k$ meeting every edge of $H$; that is, such that for every cross $((i,j),(i',j'))$, $X$ contains at least 
one of the edges $v_ju_i, u_{i'}v_{j'}$. We claim that
$G\setminus X$ is acyclic. For suppose that $C$ is a directed cycle of $G\setminus X$, with vertices $c_1\l c_t$ in order, say.
We may assume that $c_t = v_j$ say, and none of $v_1\l v_{j-1}$ are vertices of $C$. Thus $c_1\in M$, say $c_1 = u_i$.
If $c_2\in N$, say $c_2 = v_{j'}$, then $j'>j$ and so $c_2c_t\in E(G)$; but then the vertices $c_t,c_1,c_2$ are the
vertices of a directed cycle of $G$, contradicting that $G$ is $3$-free. Thus $c_2\in M$. Since $c_t\notin M$,
we may choose $s$ with $3\le s\le t$, 
minimum such that $c_s\in N$. Let $c_s = v_{j'}$, and $c_{s-1} = u_{i'}$ say. Since $c_2\l c_{s-1}\in M$ and form a directed
path in this order, and the restriction of $G$ to $M$ is acyclic, it follows that $i'>i$. Also, since
none of $v_1\l v_{j-1}$ are vertices of $C$, it follows that $j'\ge j$. If $j' = j$ then $s = t$ and $c_{t-1},c_t,c_1$
are the vertices of a directed cycle, a contradiction; so $j'>j$. Hence $((i,j),(i',j'))$ is a cross, and
$X$ contains neither of the edges $v_ju_i, u_{i'}v_{j'}$, a contradiction. This proves (1).
\bigskip

Let $(a_1,b_1)\l (a_k,b_k)$ be crosses as in (1). Let $A = \{a_1\l a_k\}$, and $B = \{b_1\l b_k\}$.
Let $C$ be the set of all $(i',j)\in V$ such that there exist $i,j'$ with $1\le i< i'\le m$ and $1\le j<j'\le n$ satisfying
$(i,j)\in A$ and $(i',j')\in B$; and let $D$ be the set of all $(i,j')\in V$ such that there exist $i',j$ with
$1\le i< i'\le m$ and $1\le j<j'\le n$ satisfying $(i,j)\in A$ and $(i',j')\in B$.
\\
\\
(2) {\em $C\cap D = \emptyset$, and $|C|+|D|\le \gamma(G)$.}
\\
\\
For suppose first that $(i,j)\in C\cap D$. Since $(i,j)\in C$, there exists $j'>j$ such that $(i,j')\in B$;
and since $(i,j)\in D$, there exists $j''<j$ such that $(i,j'')\in A$. But then $v_{j'}v_{j''}\in E(G)$
since $j''<j<j'$, and $v_{j''}u_i\in E(G)$ since $(i,j'')\in A$; and $u_iv_{j'}\in E(G)$ since
$(i,j')\in B$, contradicting that $G$ is $3$-free. This proves that $C\cap D = \emptyset$. Moreover, if
$(i',j)\in C$, we claim that $u_{i'},v_j$ are nonadjacent. For choose 
$i,j'$ with $1\le i< i'\le m$ and $1\le j<j'\le n$ such that
$(i,j)\in A$ and $(i',j')\in B$. Since $\{v_j,u_i,u_{i'}\}$ is not the vertex set of a directied cycle,
it follows that $u_{i'}v_j\notin E(G)$; and since $\{u_{i'}, v_{j'} ,v_j\}$ is also not the vertex set of a directed cycle,
$v_ju_{i'}\notin E(G)$. This proves that $u_{i'},v_j$ are nonadjacent. Similarly $u_i,v_{j'}$ are nonadjacent
for all $(i,j')\in D$. Since $C\cap D = \emptyset$, it follows that $|C|+|D|\le \gamma(G)$. This proves (2).
\bigskip

Let $a:V\rightarrow\Reals$ be defined by  $a(x) = 1$ if $x\in A$, and $a(x) = 0$ if $x\in V\setminus A$; thus,
$a$ is the characteristic function of $A$. Similarly let $b,c,d$ be the characteristic functions of $B,C,D$
respectively.
We claim that the hypotheses of \ref{ourfunctions} are satisfied. For if $1\le i<i'\le m$
and $1\le j<j'\le n$, and $a(i,j)b(i',j')>0$, then $(i,j)\in A$ and $(i',j')\in B$; hence
$v_ju_i,u_{i'}v_{j'}\in E(G)$, and so $(i',j)\in C$ and $(i,j')\in D$ from the definitions of $C,D$; and therefore
condition 1 of \ref{ourfunctions} holds. For condition 2, note first that $a(V) = k = b(V)$. Let $X,Y\subseteq V$
with $a(X) + b(Y)>a(V) = k$. We recall that $A = \{a_1\l a_k\}$ and $B = \{b_1\l b_k\}$ where $(a_i,b_i)$ is a cross
for $1\le i\le k$. Thus, $a(X) = |A\cap X|$ is the number of values of $h\in \{1\l k\}$ such that $a_h\in X$, and 
similarly $b(Y)$ is the number of $h$ with $b_h\in Y$. Since $a(X)+b(Y)>k$, there exists $h$ such that $a_h\in X$ and
$b_h\in Y$, and so $b_h$ dominates $a_h$. This proves that $b$ dominates $a$, and therefore 
the hypotheses of \ref{ourfunctions} are satisfied.

From \ref{ourfunctions}, it follows that
$a(V)b(V)\le c(V)d(V)$, and so $|A||B|\le |C||D|$. But $|A| = |B| = k$, and so $|C||D|\ge k^2$.
Consequently $|C|+|D|\ge 2k$, and hence by (2), $k\le \frac12 \gamma(G)$. This proves \ref{twocliques}.~\bbox

\section{A lemma for the second theorem}

Now we turn to the second special case of \ref{mainconj} that we can prove. The proof is in the next section, and
in this section we prove a lemma which is the main step of the proof. 
First we need some notation. Let $t\ge 1$ be an integer and let $s = 3t+1$. If $n$ is an integer, $n\mymod s$ denotes
the integer $n'$ with $0\le n'<s$ such that $n-n'$ is a multiple of $s$.
If $0\le i,j<s$ and $i,j$ are distinct, let $q>0$ be minimum such that $(i+q)\mymod s = j$ (so
$q = j-i$ if $j>i$, and $q = j-i + s$ if $j<i$). We define 
$D_s(ij) = \{(i+p)\mymod s:0\le p<q\}$.
Let $E_s$ denote the set of all ordered pairs $ij$ with $0\le i,j<s$ and $j\ne i$
such that $|D_s(ij)|\le t$, and let
$F_s$ be the set of all unordered pairs $\{i,j\}$ such that $0\le i,j<s$ and $j\ne i$ and $ij,ji\notin E_s$.
For $0\le k<s$, let $C_s(k)$ be the set of all pairs $ij\in E_s$ such that $k\in D_s(ij)$.

The lemma asserts the following.

\begin{thm}\label{circlemma}
Let $t>0$ be an integer, let $s = 3t+1$, and for $0\le i<s$ let $n_i\in \Reals$.
Then there exists $k$ with $0\le k <s$ such that 
$$\sum_{ij\in C_s(k)}n_in_j \le \frac12 \sum_{\{i,j\}\in F_s}n_in_j.$$
\end{thm} 
\Proof
Let $Q_s$ be the set of all sequences $(n_0\l n_{s-1})$ of members of $\Reals$.
We say that $(n_0\l n_{s-1})\in Q_s$ is {\em good} if there exists $k$ with $0\le k<s$ such that
$$\sum_{ij\in C_s(k)}n_in_j \le \frac12 \sum_{\{i,j\}\in F_s}n_in_j.$$
Thus we must show that every member of $Q_s$ is good. We prove this by induction on $t$.
\\
\\
(1) {\em If $t = 1$ then every member of $Q_s$ is good.}
\\
\\
For suppose that $t=1$. Let $(n_0,n_1,n_2,n_3)\in Q_s$; we must show that
there exists $k$ with $0\le k\le 3$ such that $n_kn_{k+1}\le \frac12 (n_0n_2+n_1n_3)$. 
But 
$$\min(n_0n_1,n_2n_3)^2\le n_0n_1n_2n_3 \le n_0n_1n_2n_3 + \frac14 (n_0n_2-n_1n_3)^2 = \frac14 (n_0n_2+n_1n_3)^2$$
and the claim follows. This proves (1).
\bigskip

Henceforth we assume that $t>1$.
\\
\\
(2) {\em If $(n_0\l n_{s-1})\in Q_s$ and some $n_i = 0$ then $(n_0\l n_{s-1})$ is good.}
\\
\\
For we may assume that $n_0 = 0$, from the symmetry. Define $m_i$ for $0\le i\le 3t-3$ as follows.
\begin{eqnarray*}
m_0 &=& n_{3t};\\
m_i &=& n_i\text{ for } 1\le i\le t-1;\\
m_t &=& n_t+n_{t+1};\\
m_i &=& n_{i+1}\text{ for }t+1\le i\le 2t-2;\\
m_{2t-1} &=& n_{2t}+n_{2t+1};\\
m_i &=& n_{i+2}\text{ for } 2t\le i\le 3t-3.
\end{eqnarray*}
From the inductive hypothesis and since $t>1$, the sequence $(m_0\l m_{3t-3})\in Q_{s-3}$ satisfies the theorem, and so
there exists $k'$ with $0\le k'<s-3$ such that
$$\sum_{ij\in C_{s-3}(k')}m_im_j \le \frac12 \sum_{\{i,j\}\in F_{s-3}}m_im_j.$$
If $0\le k'< t$, let $k = k'$; if $t\le k'<2t-1$, let $k = k'+1$; and if 
$2t-1\le k'\le 3t-3$, let $k = k'+2$. Since $n_0= 0$, in each case it follows easily 
(we leave checking this to the reader) that
$$\sum_{ij\in C_s(k)}n_in_j\le \sum_{ij\in C_{s-3}(k')}m_im_j.$$
But 
$$\sum_{\{i,j\}\in F_{s-3}}m_im_j = \sum_{\{i,j\}\in F_s}n_in_j - n_tn_{2t+1}\le \sum_{\{i,j\}\in F_s}n_in_j,$$
as we can check by rewriting the left side in terms of the $n_i$'s and expanding and using that $n_0 = 0$.
Consequently,
$$\sum_{ij\in C_s(k)}n_in_j\le \sum_{ij\in C_{s-3}(k')}m_im_j\le \frac12 \sum_{\{i,j\}\in F_{s-3}}m_im_j\le \frac12 
\sum_{\{i,j\}\in F_s}n_in_j,$$
and so $(n_0\l n_{s-1})$ is good. This proves (2).
\\
\\
(3) {\em Let $(n_0\l n_{s-1})\in Q_s$, such that 
$$\sum_{ij\in C_s(3t)}n_in_j\le \sum_{ij\in C_s(k)}n_in_j$$
for all $k$ with $0\le k\le 3t$. Then 
$$\sum_{0\le i< t}(t-i)(n_{3t-i} + n_{i}) \le \frac12 t\sum_{0\le i<s}n_i.$$}
For let $0\le k\le t-1$. For $0\le i\le k$, define 
$$a_i = \sum_{k+1\le j\le i+t}n_j - \sum_{i+2t+1\le j\le 3t}n_j.$$
Then 
$$\sum_{ij\in C_s(k)}n_in_j- \sum_{ij\in C_s(3t)}n_in_j = \sum_{1\le i\le k}a_in_i.$$
Since the left side of this is nonnegative, and $a_0\le a_1\le \cdots \le a_k$, it follows that
$a_k\ge 0$, that is, 
$$\sum_{k+1\le j\le k+t}n_j - \sum_{k+2t+1\le j\le 3t}n_j\ge 0.$$
Similarly, for $2t+1+k\le i\le 3t$ let
$$b_i = \sum_{i-t\le j\le 2t+k}n_j - \sum_{0\le j\le i-2t-1}n_j;$$
then 
$$\sum_{ij\in C_s(2t+k)}n_in_j- \sum_{ij\in C_s(3t)}n_in_j = \sum_{2t+1+k\le i\le 3t}b_in_i.$$
Since $b_{3t}\le b_{3t-1}\le \cdots \le b_{2t+1+k}$,
we deduce similarly that
$b_{2t+1+k}\ge 0$, that is,
$$ \sum_{t+1+k\le j\le 2t+k}n_j - \sum_{0\le j\le k}n_j\ge 0.$$
Hence
$$\sum_{k+1\le j\le k+t}n_j - \sum_{k+2t+1\le j\le 3t}n_j + \sum_{k+t+1\le j\le k+2t}n_j - \sum_{0\le j\le k}n_j\ge 0,$$
that is,
$$\sum_{k+2t+1\le j\le 3t}n_j + \sum_{0\le j\le k}n_j\le \sum_{k+1\le j\le k+2t}n_j.$$
But the sum of the left and right sides of this inequality equals $N$, 
where $N = \sum_{0\le i\le 3t}n_i$, and so the left side is at most $\frac12 N$.
Summing over all $k$ with $0\le k\le t-1$, we deduce that
$$\sum_{0\le i< t}(t-i)(n_{3t-i}+n_{i}) \le \frac12 Nt.$$
This proves (3).
\bigskip

Now to complete the proof, 
let $(n_0\l n_{s-1})\in Q_s$. 
Choose $h$ with $0\le h<s$ such that $n_h\le n_i$ for all $i$ with
$0\le i<s$. Let $n_h = x$, and for $0\le i<s$, define $m_i = n_i-x$. Thus $(m_0\l m_{s-1})\in Q_s$.
We may assume that
$$\sum_{ij\in C_s(3t)}m_im_j\le \sum_{ij\in C_s(k)}m_im_j$$
for all $k$ with $0\le k\le 3t$, by cyclically permuting $n_0\l n_{3t}$.
By (2), $(m_0\l m_{s-1})$ is good, since $m_h = 0$. Hence 
$$\sum_{ij\in C_s(3t)}m_im_j\le  \frac12 \sum_{\{i,j\}\in F_s}m_im_j.$$
But 
\begin{eqnarray*}
\sum_{ij\in C_s(3t)}n_in_j &=& \sum_{ij\in C_s(3t)}(m_i+x)(m_j+x)\\
&=& \sum_{ij\in C_s(3t)}m_im_j + \sum_{0\le k< t}x(t-k)(m_{3t-k}+m_{k}) + |C_s(3t)|x^2\\
&\le& \sum_{ij\in C_s(3t)}m_im_j + \frac12 xtM + \frac12 t(t+1)x^2,
\end{eqnarray*}
by (3), where $M = \sum_{0\le i\le 3t}m_i$. Moreover,
\begin{eqnarray*}
\frac12 \sum_{\{i,j\}\in F_s}n_in_j &=& \frac12 \sum_{\{i,j\}\in F_s}(m_i+x)(m_j +x)\\
&=& \frac12 \sum_{\{i,j\}\in F_s}m_im_j + \frac12 xtM + \frac14 stx^2\\
&\ge& \sum_{ij\in C_s(3t)}m_im_j + \frac12 xtM + \frac14 stx^2\\
&\ge& \sum_{ij\in C_s(3t)}n_in_j - (\frac12 xtM + \frac12 t(t+1)x^2) + (\frac12 xtM + \frac14 stx^2)\\
&\ge& \sum_{ij\in C_s(3t)}n_in_j.
\end{eqnarray*}
It follows that $(n_0\l n_{3t})$ is good. This completes the proof of \ref{circlemma}.~\bbox

\section{Circular interval digraphs}

We say that a digraph $G$ is a {\em circular interval digraph} if its vertices can be
arranged in a circle such that for every triple $u,v,w$ of distinct vertices, if $u,v,w$ are in clockwise order
and $uw\in E(G)$, then $uv,vw\in E(G)$. This is equivalent to saying that the vertex set of $G$ can be numbered as
$v_1\l v_n$ such that for $1\le i\le n$, the set of outneighbours of $v_i$ is $\{v_{i+1}\l v_{i+a}\}$ for some $a\ge 0$,
and the set of inneighbours of $v_i$ is $\{v_{i-b}\l v_{i-1}\}$ for some $b\ge 0$, reading subscripts modulo $n$. 
The examples given earlier to show that
conjecture \ref{mainconj} is tight infinitely often are circular interval digraphs. The main result of this section is:

\begin{thm}\label{circintthm}
$\beta(G)\le \frac12 \gamma(G)$ for every $3$-free circular interval digraph.
\end{thm}

First we need a couple of lemmas. Here is a special kind of circular interval graph. Let $t\ge 1$ be an integer,
let $n_0\l n_{3t}\ge 0$ be integers, and let $n = \sum_{0\le k\le 3t}n_i$. Let $N_0\l N_{3t}$ be disjoint sets
of cardinalities $n_0\l n_{3t}$ respectively, and let $N = N_0\cup\cdots \cup N_{3t}$. Let $N = \{v_1\l v_n\}$,
where 
$$N_i = \{v_j\;:\;n_0+n_1+\cdots + n_{i-1}< j\le n_0+n_1+\cdots + n_{i-1} +n_i\}.$$
Let $G$ be a digraph with vertex set $N$ and adjacency as follows.
\begin{itemize}
\item for $0\le k\le 3t$, if $i<j$ and $v_i,v_j\in N_k$ then $v_iv_j\in E(G)$
\item for $0\le h\le 3t$ and  $k\in \{(h+i)\mymod n\:;1\le i\le t\}$, every vertex in $N_h$ is adjacent to every vertex in
$N_k$.
\end{itemize}
In this case $G$ is a circular interval graph, and we denote it by $G(n_0\l n_{3t})$. We observe

\begin{thm}\label{bagslemma}
For all $t\ge 1$ and all choices of $n_0\l n_{3t}\ge 0$, if $G = G(n_0\l n_{3t})$ then $\beta(G)\le \frac12 \gamma(G)$.
\end{thm}
\Proof By \ref{circlemma}, there exists $k$ with $0\le k\le 3t$ such that
$$\sum_{ij\in C_s(k)}n_in_j \le \frac12 \sum_{\{i,j\}\in F_s}n_in_j,$$ 
with notation as in \ref{circlemma}. But the left side of this is at least $\beta(G)$, since every directed cycle of
$G$ contains an edge $uv$ with $u\in N_i$ and $v\in N_j$ for some $ij\in C_s(k)$; and the right side equals $\frac12 \gamma(G)$.
This proves \ref{bagslemma}.~\bbox

Let us say a $3$-free circular interval digraph is {\em maximal} if there is no pair $u,v$ of nonadjacent distinct vertices
such that adding the edge $uv$ results in a $3$-free circular interval digraph.
\begin{thm}\label{maximal}
Let $G$ be a maximal $3$-free circular interval graph. Then either $G$ is a transitive tournament, or $G$ is isomorphic
to $G(n_0\l n_{3t})$ for some choice of $t, n_0\l n_{3t}$.
\end{thm}
\Proof  Let the vertices of $G$ be $v_1\l v_n$, numbered as in the definition of a circular interval digraph, and
throughout we read these subscripts modulo $n$. For each vertex $v$, let $N^+(v), N^-(v)$ denote the set of 
outneighbours and inneighbours of $v$, respectively.
\\
\\
(1) {\em If $N^-(v)=\emptyset$ or $N^+(v)=\emptyset$ for some vertex $v$, then $G$ is a transitive tournament.}
\\
\\
For suppose that $N^-(v)=\emptyset$ for some vertex $v$, say $v_1$. If $v_kv_j\in E(G)$ for some $j,k$
with $1\le j<k\le n$, then $j>1$ and $v_1,v_j,v_k$ are in clockwise order, and therefore $v_kv_1\in E(G)$, a contradiction.
Thus $G$ is acyclic; suppose it is not a tournament. Choose $i,j$
with $1\le i<j\le n$ with $j-i$ minimum such that $v_iv_j\notin E(G)$, and let $G'$ be obtained from $G$ by adding the edge
$v_iv_j$. Then $G'$ is a $3$-free circular interval digraph, a contradiction. Thus $G$ is a tournament, and hence a 
transitive tournament since it is $3$-free. Similarly if $N^+(v)=\emptyset$ for some vertex $v$, then $G$ is a transitive
tournament. This proves (1).
\bigskip

We may therefore assume that $v_iv_{i+1}\in E(G)$ for $1\le i\le n$.
Let us say that $X\subseteq V(G)$ is a {\em cluster} if $X$ is nonempty, every two vertices in $X$ are adjacent,
$X$ can be written in the form $\{v_a,v_{a+1}\l v_b\}$ for some $a,b$,
and for every vertex $v\notin X$, either $X\subseteq N^+(v)$, or $X\subseteq N^-(v)$, or $X\cap (N^+(v)\cup N^-(v)) = \emptyset$.
\\
\\
(2) {\em For $1\le i\le n$, if $\{v_i,v_{i+1}\}$ is not a cluster, then $N^+(v_{i+1})\not \subseteq N^+(v_i)$ and
$N^-(v_{i})\not \subseteq N^-(v_{i+1})$.}
\\
\\
For certainly $v_iv_{i+1}\in E(G)$. Let $N^+(v_i) = \{v_{i+1}\l v_{i+a}\}$, where $a\ge 1$.
Suppose that $N^+(v_{i+1})\subseteq N^+(v_i)$. Then $N^+(v_{i+1}) = \{v_{i+2}\l v_{i+a}\}$.
Let the set of inneighbours of $v_i$ be $\{v_{i-b}\l v_{i-1}\}$, where $b\ge 1$, and let the set of inneighbours of
$v_{i+1}$ be $\{v_{i-c}\l v_i\}$. Thus $c\le b$; suppose that $c<b$. Then $v_{i-c-1}v_{i+1}\notin E(G)$, and also
$v_{i+1}v_{i-c-1}\notin E(G)$ since $G$ is $3$-free and $v_{i-c-1}v_i,v_iv_{i+1}\in E(G)$. Since
$v_{i-c-1}v_i,v_{i-c}v_{i+1}\in E(G)$, it follows that $v_{i-c-1}v_h, v_hv_{i+1}\in E(G)$
for all $h\in \{(i-k)\mymod n\:0\le k\le c\}$. Consequently, the digraph $G'$ obtained from
$G$ by adding the edge $v_{i-c-1}v_{i+1}$ is a circular interval digraph. From the maximality of $G$, $G'$ is not $3$-free, and
so there exists $u\in N^+(v_{i+1})\cap N^-(v_{i-c-1})$; and therefore
$u\in N^+(v_{i})\cap N^-(v_{i-c-1})$, which is impossible since $G$ is $3$-free.
This proves that $c = b$, and so $\{v_i,v_{i+1}\}$ is a cluster. Similarly if $N^-(v_{i}) \subseteq N^-(v_{i+1})$
then  $\{v_i,v_{i+1}\}$ is a cluster.
This proves (2).
\bigskip

If $X,Y$ are clusters with $X\cap Y\ne \emptyset$, it follows easily that $X\cup Y$ is a cluster. Consequently
every two maximal clusters are disjoint. Since $\{v\}$ is a cluster for every vertex $v$, it follows that the maximal clusters
form a partition of $V(G)$. Let the maximal clusters be $N_0\l N_{s-1}$ say, numbered in their natural circular order, and let
$|N_i| = n_i$
for $0\le i< s$. From the definition of a cluster, if $X,Y$ are disjoint clusters and there exists $xy\in E(G)$ with $x\in X$ and
$y\in Y$, then $xy\in E(G)$ for all $x\in X$ and $y\in Y$; we denote this by $X\rightarrow Y$. For $0\le h < s$, let 
$T_h$ be the set of all $k\in \{0\l s-1\}\setminus\{h\}$ such that $N_h\rightarrow N_k$; then $T_h = \{(h+i)\mymod s\;:1\le
i\le t_h\}$ say, for some
$t_h\ge 0$. Choose $h$ with $0\le h < s$, and choose $i$ such that $v_i\in N_h$ 
and $v_{i+1}\in N_{h+1}$. Since $\{v_i,v_{i+1}\}$ is not a cluster (because maximal clusters are disjoint), it follows from
(2) that $N^+(v_{i+1})\not \subseteq N^+(v_i)$, and so $t_{i+1}\ge t_i$. Since this holds for all choices of $i$, and
$t_0\ge t_{s-1}$, 
we deduce that
$t_0 = t_1 =\cdots = t_{s-1} = t$ say. We claim that $s = 3t+1$. For $s\ge 3t+1$ since $G$ is $3$-free; let us prove the
reverse inequality.
Let $i = n_0$ and $j = n_0 +\cdots + n_t +1$; thus $v_i\in N_0, v_{i+1}\in N_1$, $v_{j-1}\in N_t$ and $v_j\in N_{t+1}$.
Since $G$ is maximal and so adding the edge $v_iv_j$ does not result in a $3$-free circular interval digraph, 
it follows that there exists
$k$ such that $v_jv_k,v_kv_i\in E(G)$, and therefore there exists $q$ such that $q\in T_{t+1}$ and $0\in T_q$. Hence
$q-(t+1)\le t$ and $s-q\le t$; and so $s\le 3t+1$. This proves that $s  = 3t+1$, and so $G$ is isomorphic to $G(n_0\l n_{3t})$.
This proves \ref{maximal}.~\bbox

\noindent{\bf Proof of \ref{circintthm}.}
We proceed by induction on $\gamma(G)$. Suppose that $G$ is not a maximal $3$-free circular interval graph. 
Then we can add an edge to $G$
forming a $3$-free circular interval graph $G'$; and $\gamma(G')= \gamma(G)-1$, so $\beta(G')\le \frac12 \gamma(G')$
from the inductive hypothesis. Then
$$\beta(G)\le \beta(G')\le \frac12 \gamma(G')\le \frac12 \gamma(G)$$
as required.

Thus we may assume that $G$ is maximal, and we may assume that $G$ is not a transitive tournament. From \ref{maximal}
and \ref{bagslemma}, this proves \ref{circintthm}.~\bbox


\begin{thebibliography}{99}
\bibitem{AD} R. Ahlswede and D.E. Daykin, ``An inequality for the weights of two families, their unions and
intersections'',
{\em Z. Wahrscheinlichkeitsth. verw. Gebiete}, 43 (1978),  183-–185.
\bibitem{CH} L. Caccetta and R. H\"aggkvist, ``On minimal digraphs with given girth'', {\em Proceedings of the Ninth Southeastern
Conference on Combinatorics, Graph Theory and Computing} (Florida Atlantic University, Boca Raton, Florida, 1978),
Congressus Numerantium XXI, Utilitas Math., 1978, 181-187.
\bibitem{KS} A. Kostochka and M. Stiebitz, in a lecture by Kostochka at the American Institute of Mathematics in January 2006.

\end{thebibliography}
\end{document}